\documentclass[a4paper,12pt]{article}

\usepackage{graphicx}
\usepackage{amsmath} 
\usepackage{amssymb}  
\usepackage{bm} 
\usepackage{float}
\usepackage{mathtools}
\usepackage{amsthm} 
\theoremstyle{definition}
\newtheorem{example}{Example} 

\newcommand{\dif}{\mathrm{d}} 
\newcommand{\R}{\mathbb{R}}
\newcommand{\N}{\mathbb{N}}

\def\N{\mathbb{N}}
\def\R{\mathbb{R}}

\title{\bf Measures and LMI for impulsive optimal control with applications to space rendezvous problems}

\begin{document}

\author{Mathieu Claeys$^{1,2}$ Denis Arzelier$^{1,2}$\\
Didier Henrion$^{1,2,3}$ Jean-Bernard Lasserre$^{1,2,4}$}

\footnotetext[1]{CNRS; LAAS; 7 avenue du colonel Roche, F-31077 Toulouse; France.}
\footnotetext[2]{Universit\'e de Toulouse; UPS, INSA, INP, ISAE; UT1, UTM, LAAS; F-31077 Toulouse; France}
\footnotetext[3]{Faculty of Electrical Engineering, Czech Technical University in Prague,
Technick\'a 2, CZ-16626 Prague, Czech Republic}
\footnotetext[4]{Institut de Math\'ematiques de Toulouse, Universit\'e de Toulouse; UPS; F-31062 Toulouse, France.}

\maketitle

\begin{abstract}
This paper shows how to find lower bounds on, and sometimes solve globally, a large class of nonlinear optimal control problems with
impulsive controls using semi-definite programming (SDP).
This is done by relaxing an optimal control problem into a measure differential problem.
The manipulation of the measures by their moments reduces the problem to a convergent series of standard linear matrix inequality (LMI) relaxations.
After providing numerous academic examples, we apply the method to the impulsive rendezvous of two orbiting spacecrafts.
As the method provides lower bounds on the global infimum, global optimality of the solutions can be guaranteed numerically by a posteriori simulations, and we can recover simultaneously the optimal impulse time and amplitudes by simple linear algebra.
\end{abstract}

\section{Introduction}

Optimal control problems are still an active area of current research despite the availability of powerful theoretical tools such as Pontryagin's maximum principle or the Hamilton-Jacobi-Bellman approach, that both provide conditions for optimality. However, numerical methods based on such optimality conditions rely on a certain number of assumptions that are often not met in practice. In addition, state constraints are particularly hard to handle for most of the methods.

On the other side, many numerical methods have been developed that deliver locally optimal solutions satisfying sufficient optimality conditions. However, the users of these methods are often left to wonder if a better solution exists. For example, in the particular case of impulsive controls, it is often not known if more regular solutions could provide a better cost.
For a recent survey on impulsive control see e.g. \cite{Kurzhanski} and the references therein. See also \cite{rapaport} for a recent application and for more references\footnote{We are grateful to T\'erence Bayen for pointing out this reference to us.}. For historical works see e.g. \cite{Neustadt,Rischel,Schmaedeke} and also \cite{Bensoussan}.

This paper presents a method based on \cite{LasserreNLOCP,pocp} but covering a larger class of solutions, including impulsive controls. This algorithm provides a sequence of non-decreasing lower bounds on the global minimizers of affine-in-the-control polynomial optimal control problems. In particular, it may assert the global optimality of local solutions found by other methods. As importantly, the algorithm is also able to provide numerical certificates of infeasibility or unboundedness for ill-posed problems. Finally, in some cases, it is also possible to generate the globally optimal control law.

At the end of the paper, this method is successfully applied to the problem of coplanar space fuel-optimal linearised rendezvous. We show with two different examples from the literature that the proposed algorithm is able to retrieve the impulsive optimal solution conjectured by running a direct approach based on the solution of a Linear Programming (LP) problem. Without assuming the nature of the propulsion (continuous or impulsive), the obtained impulsive solution is certified to be a global fuel-optimal solution.

\subsection{Contributions}
The paper improves the method presented in \cite{LasserreNLOCP,pocp} in the following ways:
\begin{itemize}
  \item Impulsive control can now be taken into account.
  \item Because controls are represented by measures and not by variables, the size of SDP blocks is significantly reduced. This allows to handle larger problems in terms of number of state variables as well as to reach higher LMI relaxations.  
  \item Total variation constraints can be handled very easily.
\end{itemize}

These three improvements make it altogether possible to tackle problems such as consumption minimization for space rendezvous, the other significant contribution of this paper.

\subsection{Notations}

Integration of a function $f:\R^n\to\R$ with respect to a measure $\mu$ on a set $X\subset\R^n$ is written $\int_{X} \! f(x) \, \dif \mu(x)$. The Lebesgue or uniform measure on $X$ is denoted by $\lambda$ whereas
the Dirac measure concentrated at point $x$ is denoted by $\delta_x$. A measure $\mu$ is a probability measure whenever $\int d\mu = 1$.
The support of measure $\mu$ is denoted by ${\rm supp}\,\mu$. The indicator function of set $X$ (equal to one
in $X$ and zero outside) is denoted by $I_X$.

$F(X)$ is the space of Borel measurable functions on $X$, whereas $BV(X)$ is the space of functions of bounded variation on $X$.
 $\R[z]$ is the ring of polynomials in the variable $z$. $\mathcal{B}(X)$ denotes the Borel $\sigma$-algebra associated with $X$.

If $k \in \N^n$ denotes a vector of indices then $x^k$ with $x \in \R^n$ is the multi-index notation for $\prod x_i^{k_i}$. The degree of the index $k$ is $\mathrm{deg}\,k=\sum k_i$. Finally,
$\N_d^n$ is the set of all indices for which $\mathrm{deg}\,k \leq d, \; k \in \N^n$.

\section{The optimal control problem}

This paper deals with the following nonlinear optimal control problem
\begin{equation}
\begin{array}{rcl}
V(x_0) & = \displaystyle \inf_{u(t) \in F([0,T])^m} & I(x_0,u) = \displaystyle \int_0^T h(t,x(t))dt + \int_0^T H(t) u(t)dt + \,\,h_T(x(T)) \\[1em]
  & \mathrm{s.t.} &  \dot{x}(t) = f(t,x(t)) + G(t) u(t), \quad t \in [0,T] \\[.5em]
  && x(0)=x_0 \in X_0, \quad x(T) \in X_T, \quad x(t) \in X \subset \R^n.
    \label{eq:OCPo}
\end{array}
\end{equation}  
where the dot denotes differentiation w.r.t. time and the prime denotes transposition.
Criterion $I(x_0,u)$ is {\it affine} in the control $u$, and $V$ is called the value function. 
It is assumed that all problem data are polynomials, meaning that all functions are in $\R[t,x]$, and that all sets are compact basic semialgebraic.
Recall that such sets are those which may be written as
$\left\lbrace  x : \: a_i(x) \geq 0, \: i =1, \ldots ,m \right\rbrace$ for some family $\left\lbrace a_i  \right\rbrace_{j=1}^m \subset \R[x]$. A mild technical condition (implying compactness of $X$) must be satisfied \cite[Assumption 2.1]{Lasserre},
but it is often met in practice (for instance, an additional standard ball constraint $\sum x_i^2 \leq r^2$ enforces the condition).
The reason for making these assumptions will be apparent in the later sections.

Without additional assumptions and constraints, the infimum in problem (\ref{eq:OCPo})
is generally not attained in the space of measurable functions \cite{Young}. For this reason, in this paper we consider problems
for which controls are allowed to be generalized functions, i.e. {\it measures}, thereby extending the original formulation
as follows:
\begin{equation}
\begin{array}{rcl}
V_R(x_0) & = \displaystyle \inf_{w(t) \in BV([0,T])^m} & I(x_0,w) =\displaystyle \int_0^T h(t,x(t))dt + \int_0^T H(t) dw(t) + \,\,h_T(x(T)) \\[1em]
  & \mathrm{s.t.} &  dx(t) = f(t,x(t))\,dt + G(t) dw(t), \quad t \in [0,T] \\[.5em]
  && x(0)=x_0 \in X_0, \quad x(T) \in X_T, \quad x(t) \in X \subset \R^n
    \label{eq:OCP}
\end{array}
\end{equation}  
where $V_R$ stands for the relaxed value function.
In particular, in problem (\ref{eq:OCP}) controls may be impulsive: the (vector) control can been seen as a (vector) distribution
of the first order and it is therefore the distributional derivative $dw(t)$ of some (vector) function of bounded
variation $w(t) \in BV([0,T])^m$, see e.g. \cite{Schmaedeke}
and \cite[\S 4]{Riesz} or also \cite[Prop. 8.3]{Brezis}.

\section{The measure problem}
\label{sec:measureProblem}

In this section, we formulate problem (\ref{eq:OCP}) into a measure differential problem, a necessary step towards obtaining a tractable SDP problem. 
Optimal control problems involving measures have been introduced to accept solutions that are ruled out or ill-defined
in classical optimal control, see e.g. \cite{Young}. Multiple solutions, impulsive or chattering controls can be handled naturally by the associated measure problem. This section, rather than providing rigorous proofs, outlines the main ideas behind this transformation.

A few remarks are worth pointing out. First of all, it is crucial that $G(t)$ be a matrix of smooth functions, an hypothesis automatically fulfilled by polynomials. As a matter of fact, multiplying distributions with such functions is a well-defined operation (unlike e.g. the product of two distributions). Therefore, except for some very particular cases \cite{Miller}, $G$ cannot be a function of states $x_j$ that could potentially present jump discontinuities. To simplify notations, we have simply assumed that $G$ depends on $t$ only\footnote{In all rigour, it could be possible to include state jumps in $G$, but this requires a careful definition of what is meant by integration, as done e.g. for studying stochastic differential equations. This goes well beyond the scope of this paper.}.
Secondly, in the absence of impulses, the distributional differential is the traditional differential, and the dynamics are classical differential equations with controls
$dw(t)=u(t)dt$ which are absolutely continuous with respect to the Lebesgue measure. Finally, it must be noted that state trajectories $x(t)$ are themselves functions of bounded variations, being the sum of two such functions, and that this is their broadest class in the sense that there is no more general distribution solutions for the states \cite{Schmaedeke}.

Because distributional derivatives of functions of bounded variation on compact supports can be identified with measures \cite[\S 50]{Riesz}, the dynamics in problem (\ref{eq:OCP}) may be interpreted as a measure differential equation. As $X\subset\R^n$ is assumed to be compact, by one of the Riesz representation theorems \cite[\S 36.6]{Kolmogorov}, these measures can be put in duality correspondence with all continuous functions $v(t,x(t))$ supported on $[0,T] \times X$. We will use these test functions to define linear relations between the measures. Note that because continuous functions on compact sets can be uniformly approached by polynomials by virtue of the Stone-Weierstrass theorem, it is enough to consider polynomial test functions $v(t,x(t)) \in \R [[0,T] \times X]$.

By Lebesgue's decomposition theorem \cite[\S 33.3]{Kolmogorov}, we can split the control measures $w(dt)$ into two parts: their absolutely continuous parts with density $u:[0,T]\to\R^m$ (with respect to the Lebesgue measure) and their purely singular parts with {\it jump amplitude} vectors $u_{t_j}\in\R^m$ supported at impulsive {\it jump instants} $t_j$, $j\in J$, with $J$ a subset of Lebesgue measure zero of $[0,T]$, not necessarily countable\footnote{We suspect however that for the control problems studied in this paper, subset $J$ can be assumed countable without loss of generality.}. We write
\[
w(dt) = u(t)dt + \sum_{j \in J} G(t_j) u_{t_j} \delta_{t_j}(dt)
\]
to model jumps in state-space
\[
x^+(t_j) = x^-(t_j) + G(t_j) u_{t_j},\quad \forall\,j\in J.
\]

Now, {\it given an initial state} $x_0\in X_0$ and {\it given a control} $w(t)\in BV([0,T])^m$, denote by $x(t) \in BV([0,T])^n$ the corresponding feasible trajectory.
Then for smooth test functions $v:[0,T]\times\R^n\to\R$, it holds
\begin{equation}\label{time}
\begin{array}{rcl}
  \displaystyle \int_0^T dv(t,x(t)) & = & v(T,x(T)) - v(0,x(0))\\
  & = & \displaystyle \int_0^T \!\! \left( \frac{\partial v}{\partial t}+\left(\frac{\partial v}{\partial x}\right)' f \right) dt\\
  && +\,\displaystyle \int_0^T \! \left(\frac{\partial v}{\partial x}\right)' G u dt\\
  && +\,\displaystyle \sum_{j\in J}v(t_j,x^+(t_j))-v(t_j,x^-(t_j)).
\end{array}
\end{equation}
We are going to express the above {\it temporal} integration (\ref{time}) along the trajectory
in terms of {\it spatial} integration with respect to appropriate and so-called {\it occupation measures}.
For this purpose, define:
\begin{itemize}
\item The {\it time-state} occupation measure
  \[
   \mu[x_0,w(t)](A\times B) = \int_A I_B(x(t))\,dt,\quad
   \forall A\in\mathcal{B}([0,T]), \quad \forall B \in \mathcal{B}(X)
  \]
which measures the {\it occupation} of $A\times B$ by the pair $(t,x(t))$ all along the trajectory.
Note that we write $\mu[x_0,w(t)]$ to emphasize the dependence of $\mu$ on initial condition $x_0$ and
control $w(t)$. However, for notational simplicity, we may use the notation $\mu$.
By a standard result on Borel measures on a cartesian product, the occupation measure $\mu$ can be {\it disintegrated} into
\[
\mu(A\times B) = \int_{A}\xi(B\,\vert\,t)\,dt,
\]
where $\xi(dx\,\vert\,t)$ is the distribution of $x\in\R^n$, conditional on $t\in [0,T]$. It is a stochastic kernel, i.e.,
\begin{itemize}
\item for every $t\in [0,T]$, $\xi(\cdot\,\vert\,t)$ is a probability distribution on $X$, and
\item for every $B\in \mathcal{B}(X)$, $\xi(B\,\vert\,\cdot)$ is a Borel measurable function on $[0,T]$.
\end{itemize}
In our case, since the initial state $x_0$ and the control $w(t)$ are given,
the stochastic kernel $\xi(dx\,\vert\,t)$ is well defined along continuous
arcs of the trajectory as
\begin{equation}\label{eq:dirac}
\xi(B|t) = I_B(x(t)) = \delta_{x(t)}(B),\quad\forall t \in [0,T]\setminus J,\quad\forall B\in\mathcal{B}(X).
\end{equation}  
On the other hand, at every jump instant $t_j\in J$, we let
\[
\xi(B\,\vert\,t_j)=\frac{\lambda(B\cap[x^-(t_j),x^+(t_j)])}{\lambda([x^-(t_j),x^+(t_j)])},
\quad \forall t_j\in J, \quad\forall B\in\mathcal{B}(X).
\]
This means that the state is uniformly distributed along the segment
linking the state before and after the jump, the above denominator ensuring
that $\xi(\cdot\,\vert\,t)$ has unit mass for all $t$.
\item The {\it control-state} occupation measure
\begin{equation*}
  \nu[x_0,w(t)](A \times B)\,=\,\int_A \xi(B\,\vert\,t)\,dw(t),\quad
\forall A\in\mathcal{B}([0,T]), \quad \forall B \in \mathcal{B}(X).
\end{equation*}
\item The {\it final state} occupation measure
\[
 \mu_T[x_0,w(t)](B) = I_B(x(T)),\quad \forall B\in\mathcal{B}(X_T).
\]
\end{itemize}
With these definitions, Eq. (\ref{time}) may be written in terms of measures as:
\begin{align}
  \int_{X_T} v(T,x) \, d\mu_T(x) - v(0,x_0)  = & \label{eq:occup} \\
  \int_{[0,T] \times X} \left(\frac{\partial v}{\partial t} + \left(\frac{\partial v}{\partial x}\right)' f\right)\,d\mu(t,x) \: + &
  \int_{[0,T] \times X} \left(\frac{\partial v}{\partial x}\right)' G d\nu(t,x) = \notag \\
  \int_{[0,T]}\left[\int_X \left(\frac{\partial v}{\partial t} + \left(\frac{\partial v}{\partial x}\right)'\left(f +
  G u\right)\right)\,\xi(dx\,\vert\,t)\right]dt \: + &
  \sum_{j\in J}v(t_j,x^+(t_j))-v(t_j,x^-(t_j)). \notag
\end{align}
Similarly, the criterion in (\ref{eq:OCP}) to evaluate the trajectory and the control reads:
\[
  I(\mu,\nu,\mu_T) = \int_{[0,T] \times X} \! h \, d\mu +
              \int_{[0,T] \times X} \! H \, d\nu +
              \int_{X_T} \! h_T \, d\mu_T.
\]
In view of the above formulation with occupation measures, one may now define a {\it relaxed version} (or {\it weak} formulation)
of the initial (measure) control problem (\ref{eq:OCP}). First note that
\[
V_R(x_0)=\inf_{w(t)} I(\mu[x_0,w(t)],\nu[x_0,w(t)],\mu[x_0,w(t)])
\]
where the infimum is taken over all the occupation measures defined above, 
corresponding to a given initial condition $x_0$ and control $w(t)$.
Second, instead of searching for a control $w(t)$, we search for 
a triplet of measures that solves the infinite dimensional problem:
\[
 V_M(x_0) = \inf_{\mu,\nu,\mu_T} I(\mu[x_0],\nu[x_0],\mu_T[x_0])
\]
under the {\it trajectory} constraints (\ref{eq:occup}) for all
$v \in \R[t,x]$ and the {\it support} constraints
${\rm supp}\,\mu\,=\,{\rm supp}\,\nu\,=\,[0,T]\times X$, ${\rm supp}\,\mu_T=X_T$.
The measures now depend only on initial condition $x_0$, since they just
have to satisfy linear constraints (\ref{eq:occup}). This motivates
the notation $\mu[x_0]$, $\nu[x_0]$, $\mu_T[x_0]$ in the above problem.
This problem is an obvious relaxation of problem (\ref{eq:OCP}) 
which is itself a relaxation of (\ref{eq:OCPo}), hence
\[
 V_M(x_0)\leq V_R(x_0)\leq V(x_0).
\]
In the remainder of the paper, we will deal with this relaxed version of
the occupation measures problem.
However, for a well-defined control problem (\ref{eq:OCP}) one expects that in fact
$V_M(x_0)= V_R(x_0)$ and that an optimal solution of the relaxed problem
will be the triplet of occupation measures corresponding to an optimal trajectory of problem (\ref{eq:OCP})
with given initial state $x_0$ and control $w(t)$.
Note that for the standard polynomial optimal control problem (\ref{eq:OCPo}), without impulsive
controls, and under additional convexity assumptions, it has been proved in \cite{LasserreNLOCP}
that indeed $V_M(x_0)= V_R(x_0) = V(x_0)$.

\subsection{Initial state with a given distribution}

Recall that the occupation measures defined in the previous section all depend on $x_0$.
Observe that if $\mu_0$ is a given probability measure on $X_0\subset\R^n$ and if one now defines:
\[
\begin{array}{rcl}
\mu(A\times B) & = & \int_{X_0}\mu[x_0](A\times B)\,d\mu_0(x_0), \\
\nu(A\times B) & = & \int_{X_0}\nu[x_0](A\times B)\,d\mu_0(x_0), \\
\mu_T(B) & = & \int_{X_0}\mu_T[x_0](B)\,d\mu_0(x_0)
\end{array}
\]
for all $A\in\mathcal{B}([0,T])$ and $B \in \mathcal{B}(X)$, then
\[
  I(\mu[\mu_0], \nu[\mu_0], \mu_T[\mu_0]) = 
  \int_{X_0} I(\mu[x_0],\nu[x_0],\mu_T[x_0])d\mu_0(x_0)
\]                           
becomes the expected average cost associated 
with the trajectories and with respect to the probability measure $\mu_0$ on $X_0$.

Therefore, the relaxed problem with measures now reads as follows:
\begin{equation}\label{eq:OCPM}
\begin{array}{rc@{}l}
 V_M (\mu_0) = & \displaystyle \inf_{\mu,\nu,\mu_T} & I(\mu[\mu_0], \nu[\mu_0], \mu_T[\mu_0]) = \displaystyle \int h d\mu + \int H d\nu + \int h_T d\mu_T \\[1em]
  & \mathrm{s.t.} & \displaystyle   \int v \, d\mu_T - \int v \, d\mu_0 = \int \left(\frac{\partial v}{\partial t} + \left(\frac{\partial v}{\partial x}\right)' f\right)\,d\mu
 + \int \left(\frac{\partial v}{\partial x}\right)' G d\nu \\[1em]
  && {\rm supp}\,\mu\,=\,{\rm supp}\,\nu\,=\,[0,T]\times X,\quad {\rm supp}\,\mu_T=X_T.
\end{array}
\end{equation}
Note that in this case, the stochastic kernel $\xi(dx|t)$ along continuous arcs
of the trajectory is generally not a Dirac measure as in (\ref{eq:dirac}),
unless $\mu_0$ is a Dirac measure supported at $x_0$ and
the optimal control $w$ is unique.

By solving this relaxed problem we expect that its optimal value satisfies
\[
V_M(\mu_0) = \int_{X_0} V_M(x_0)d\mu_0(x_0),
\]
i.e. that $V_M(\mu_0)$ is
the expected average cost associated 
with optimal trajectories and with respect to the probability measure $\mu_0$ on $X_0$.

\subsection{Free initial state}

In this case, in addition to the control we also have the freedom of choosing the 
best possible initial state. For this purpose introduce an {\it unknown}  probability measure $\mu_0$ on $X_0$. 
Then the relaxed problem with measures, analogue of (\ref{eq:OCPM}), 
reads almost the same except that:
\begin{itemize}
\item we now optimize over $\mu,\nu,\mu_0,\mu_T$;
\item in the support constraints we introduce the additional constraint
${\rm supp}\,\mu_0=X_0$.
\end{itemize}
By solving this relaxed problem we now expect that its optimal value
denoted $V_M(X_0)$ satisfies:
\[ 
V_M(X_0) = \inf_{\mu_0} V_M(\mu_0) = \inf_{x_0 \in X_0} V_M(x_0).
\]

\subsection{Decomposition of control measures}
\label{sec:decomposition}
All measures in (\ref{eq:OCPM}) are positive measures, except for the signed measures $\nu$ which deserve special treatment for our purposes.
Using the Jordan decomposition theorem \cite[\S 34]{Kolmogorov}, these measures may be split into a positive part $\nu^+$ and negative part $\nu^-$, that is $\nu = \nu^+ - \nu^-$, both being positive measures.

This decomposition has the added benefit of providing an easy expression for the $L_1$ norm of the control, which is sometimes to be constrained or optimized in some problems. Indeed, define the {\it total variation} control measure by
\begin{equation*}
  | \nu | = \nu^+ + \nu^-.
\end{equation*}
The total variation norm of the measure $\nu$ is just the mass of $\vert\nu\vert$, i.e.,
\begin{equation*}
  \|\nu\|_{TV} = \int \! d | \nu |.
  \label{eq:totalVariation}
\end{equation*}

\subsection{Handling discrete control sets}
\label{sec:discreteControl}
It is often desirable to restrict the set of admissible controls to be a subset of $\R$. Here we will limit ourselves to the very important case of handling discrete control sets. Let us assume that controls $u$ are only allowed to take their values in $U=\left\lbrace u_1,..., u_m \right\rbrace$. Define $\nu_i$ as the probability measures of choosing controls $u_i$. Then clearly, the total probability of choosing one of the controls in $U$ must be 1 at each time along the trajectory.
Then the control measures $\nu$ are simply the linear combination of the probability measures weighted by their respective control values:
\begin{equation*}
  \nu = \sum_i u_i \, \nu_i.
\end{equation*} 
Using the same method as in \S\ref{sec:measureProblem}, we have $\forall \, v(t) \in \R[t]$:
\begin{equation*}
  \sum_i \int \!  v(t) \, u_i \, d\nu_i(t,x) = \int \! v(t) \, d\nu(t,x).
\end{equation*}

Note that with this substitution, all measures involved in the measure problem are now positive; there is no need to implement the trick of \S\ref{sec:decomposition}. Using these extra constraints, it is now possible to solve bang-bang control problems.

\subsection{Summary}

To summarize, the advantages for introducing the relaxed control problem (\ref{eq:OCPM})
with measures are the following:
\begin{itemize}
\item controls are allowed to be measures with absolutely continuous components and singular components including impulses;
\item state constraints are easily handled via support constraints;
\item the initial state has a fixed given distribution on some pre-specified domain;
\item a free initial state in some pre-specified domain is also allowed.
\end{itemize}

\section{The moment problem}
So far, the hypothesis of polynomial data has not been used, but its crucial importance will appear in this section, where measures will be manipulated through their moments. This will lead to a semi-definite programming (SDP) problem featuring countably many equations.

Define the moments of measure $\mu$ as
\begin{equation}\label{moments}
y_k^\mu = \int_X \! z^k \, d\mu(z).
\end{equation}
Then, with a sequence $y=(y_k)$, $k \in \N^n$, let $L_y:\R[z]\to\R$ be the linear functional
\begin{equation*}
  f \; \left( = \sum_k f_k z^k \right)
  \quad \mapsto \quad
  L_y(f)\,=\,\sum_k f_k y_k,\quad f\in\R[z].
\end{equation*}
Define the moment matrix of order $d \in \N$ associated with $y$
as the real symmetric matrix $M_d\left(y\right)$ whose $(i,j)$th entry reads
\begin{equation*}
  M_d(y)[i,j] = L_y \left(z^{i+j}\right) = y_{i+j},
  \quad \forall i, j \in \N^n_d.
\end{equation*}
Similarly, define the localizing matrix of order $d$ associated with $y$ and $h \in \R[z]$ as the real symmetric matrix $M_d(h\,y)$ whose $(i,j)$th entry reads
\begin{equation*}
  M_d(h \, y)[i, j] = L_y \left( h(z) \, z^{i+j} \right) = \sum_k h_k \, y_{i+j+k},
  \quad \forall i,j \in \N^n_d.
\end{equation*}
As a last definition, a sequence $y^\mu=(y_{k}^\mu)$ is said to have a representing measure if there exists a finite Borel measure $\mu$ on $X$, such that
relation (\ref{moments}) holds for every $k \in \N^n$.

Now comes the crucial result of the section: a sequence of moments $y^\mu$ has a representing measure defined on a semi-algebraic set $X^\mu = \{x \: :\:
a^\mu_i(x) \geq 0, \: i=1,2,\ldots\}$ if and only if $M_d(y^\mu) \succeq 0,\; \forall \, d \in \N$ and $M_d(a^\mu_i \, y^\mu) \succeq 0, \; \forall \, d \in \N$ and $\forall a^\mu_i$ defining set $X^\mu$ \cite[Theorem 3.8]{Lasserre}. This has the very practical implication that the measure problem defined in (\ref{eq:OCPM}) has an equivalent formulation in terms of moments. Indeed, because all problem data were assumed to be polynomial, the criterion in (\ref{eq:OCPM}) can be transformed into a linear combination of moments to be minimized:
\begin{equation}
  V_m = \inf_y \,\, (b^\mu)' y^\mu
        + (b^\nu)' y^\nu
        + (b^{\mu_T})' y^{\mu_T} = b'y
  \label{eq:momentCriterion}
\end{equation}
where the infimum is now over the aggregated sequence $y$ of moments of all the measures.
Because the test functions were also restricted to be polynomials, the constraints in (\ref{eq:OCPM}) can be turned into countably many linear constraints on the moments: 
\begin{equation}
  A^\mu y^\mu + A^\nu y^\nu + A^{\mu_0} y^{\mu_0} + A^{\mu_T} y^{\mu_T} = Ay = 0.
  \label{eq:momentConstraints}
\end{equation}
The only non-linear part are the SDP constraints for measure representativeness, to be satisfied $\forall d \in \N$:
\begin{align}
  & M_d(y^\mu) \succeq 0, \quad M_d(a_i^\mu \, y^\mu) \succeq 0, \notag \\
  & M_d(y^\nu) \succeq 0, \quad M_d(a_i^\nu \, y^\nu) \succeq 0, \notag \\ \label{eq:momentRepr} 
  & M_d(y^{\mu_0}) \succeq 0, \quad M_d(a_i^{\mu_0} \, y^{\mu_0} ) \succeq 0, \\
  & M_d(y^{\mu_T}) \succeq 0, \quad M_d(a_i^{\mu_T} \, y^{\mu_T} ) \succeq 0. \notag
\end{align}

\section{LMI relaxations}
The final step to reach a tractable problem is relatively obvious: we simply truncate the problem to its first few moments. Let $d_1 \in \N$ be the smallest integer such that all criterion monomials belong to $\N^{n+1}_{2 d_1}$. This is the degree of the so called \emph{first relaxation}. For each relaxation, we reach a standard LMI problem that can be solved numerically by off-the-shelf software by simply truncating Eq. (\ref{eq:momentCriterion}), (\ref{eq:momentConstraints}) and (\ref{eq:momentRepr}) to involve only moments in $\N^{n+1}_{2 d}$, with $d \geq d_1$ the relaxation order.

Observe that $d_j > d_i \; \Rightarrow V_M^{d_j} \geq V_M^{d_i} $. Therefore, by solving the truncated problem for ever greater relaxation orders, we will obtain a monotonically non-decreasing sequence of lower bounds to the true cost. In the examples below, we will see that in practice, the optimal cost is usually reached after a few relaxations.

\section{Academic examples}
In this section, many examples are presented to showcase the different features of the method.
Ex. \ref{ex:basic} to \ref{ex:unbounded} are variations of the same basic problem to give a thorough tour of the method's capabilities. Ex. \ref{ex:VDP} is taken from the literature and shows how the method compares to, or rather nicely complements, existing optimal control algorithms. All examples use GloptiPoly \cite{GloptiPoly} for building the truncated LMI moment problems and SeDuMi \cite{SeDuMi} for their numerical solution.

Before proceeding to the examples, define the marginal $M_d \left( y, z \right)$ of a moment matrix with respect to variable $z$ as the moment matrix of the subsequence of moments concerning polynomials of $z$ only.

\begin{example}[Basic impulsive problem]
\label{ex:basic}

\begin{equation*}
  V = \inf_{u(t)} \! \int_0^2 x^2(t) \, \dif t
\end{equation*}
such that
\begin{align*}
    & \dot{x}(t) = u(t) \\
    & x(0) = 1, \quad x(2) = \frac{1}{2} \\
    & x^2(t) \leq 1.
\end{align*}

In this introductory example, it is straightforward to notice that the optimal solution consists of
reaching the turnpike $x(t)=0$ by an impulse at initial time $t=0$, and likewise, departing from it by an impulse at final time $t=T=2$, see Fig. \ref{fig:basic}.

\begin{figure}
  \centering
  \includegraphics[width=0.7\textwidth]{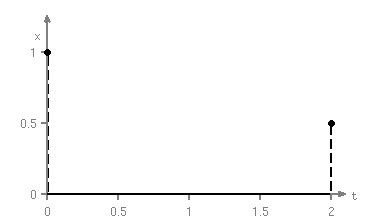}
  \caption{Trajectory for Ex. \ref{ex:basic} }
  \label{fig:basic}
\end{figure}

The associated measure problem reads:
\begin{equation*}
  V_M = \inf_{\mu, \nu} \int_{[0,T] \times X} \! x^2 \: \dif \mu
\end{equation*}
such that
\begin{align*}
    & \int_{X_T} \! v \, \dif \mu_T - \int_{X_0} \! v \, \dif \mu_0 =
           \int_{[0,T] \times X} \! \frac{\partial v}{\partial t} \, \dif \mu +
           \int_{[0,T] \times X} \! \frac{\partial v}{\partial x} \, \dif \nu \quad \forall \, v \in \R[t,x] \\
    & \mu_0 = \delta_0 \quad X_0 = \{1\} \quad \mu_T  = \delta_{\frac{1}{2}} \quad X_T = \left\{\frac{1}{2}\right\} \\
    & X = \left\lbrace x \in \mathbb{R}:  1 - x^2 \geq 0 \right\rbrace.
\end{align*}

Using the procedure outlined above, one obtains a series of truncated moment problems that can be solved by semi-definite programming.
Letting $y_{ij}^\mu=\int \! t^i x^j  \, d\mu$, the first LMI relaxation is 
\begin{equation*}
  V_M^1 = \inf_y y_{02}^\mu
\end{equation*}
subject to the linear constraints associated to the dynamics:
\begin{align*}
     y_{00}^{\mu_T} - y_{00}^{\mu_0} & = 0 \\
     y_{10}^{\mu_T} - y_{10}^{\mu_0} & = y_{00}^\mu \\
     y_{01}^{\mu_T} - y_{01}^{\mu_0} & = y_{00}^{\nu^+} - y_{00}^{\nu^-} \\
     y_{20}^{\mu_T} - y_{20}^{\mu_0} & = 2 y_{10}^\mu \\
     y_{11}^{\mu_T} - y_{11}^{\mu_0} & = y_{01}^\mu + y_{10}^{\nu^+} - y_{10}^{\nu^-} \\
     y_{02}^{\mu_T} - y_{02}^{\mu_0} & = 2 y_{01}^{\nu^+} - 2 y_{01}^{\nu^-},    
\end{align*}
to the SDP representativeness constraints for $\tau = \left\lbrace \mu, \nu^+, \nu^- \right\rbrace$:
\begin{equation*}
     \begin{bmatrix}
       y_{00}^\tau & y_{10}^\tau & y_{01}^\tau \\
       y_{10}^\tau & y_{20}^\tau & y_{11}^\tau \\
       y_{01}^\tau & y_{11}^\tau & y_{02}^\tau \\
     \end{bmatrix} \succeq 0, \qquad y_{00}^\tau - y_{02}^\tau \geq 0,
\end{equation*}
and to the boundary conditions:
\begin{align*}
 \begin{bmatrix}
       y_{00}^{\mu_0} & y_{10}^{\mu_0} & y_{01}^{\mu_0} & y_{20}^{\mu_0} & y_{11}^{\mu_0} & y_{02}^{\mu_0}
  \end{bmatrix}
     = 
     \begin{bmatrix}
       1 & 0 & 1 & 0 & 0 & 1
     \end{bmatrix}, \\
     \begin{bmatrix}
     y_{00}^{\mu_T} & y_{10}^{\mu_T} & y_{01}^{\mu_T} & y_{20}^{\mu_T} & y_{11}^{\mu_T} & y_{02}^{\mu_T}
     \end{bmatrix}
     = 
     \begin{bmatrix}
       1 & 2 & \frac{1}{2} & 4 & 1 & \frac{1}{4}
     \end{bmatrix}. &
\end{align*}

It turns out that the optimal value $V_M=0$ is estimated correctly (to numerical tolerance) from the first relaxation on and that the optimal trajectory $x(t)=0$ can easily be recovered. Indeed, the marginal $M_d \left( y^\mu, x  \right)$ is the length of the time interval multiplying a truncated moment matrix of a Dirac measure concentrated at $x=0$, while its marginal with respect to $t$ equals a truncated Lebesgue moment matrix on the $[0,2]$ interval. More importantly, one can recover the optimal controls as the marginal $M_d(y^\nu,t)$ is the weighted sum of Dirac measures located at the impulse times, the weights being the impulse amplitudes. In summary, we can recover numerically the optimal measures
\[
\mu(dt,dx)=I_{[0,2]}(dt)\delta_0(dx), \quad
\nu(dt,dx)=-\delta_0(dt)I_{[0,1]}(dx)+\delta_2(dt)I_{[0,\frac{1}{2}]}(dx).
\]
\end{example}

\begin{example}[Total variation constraints]
\label{ex:totalVariation}
We take back Ex. \ref{ex:basic} with an additional constraint on the total variation of the control:
\begin{equation*}
  \int_0^2 \! |u(t)| \, \dif t  \leq 1
\end{equation*}
whose measure equivalent reads:
\begin{equation*}
 \|\nu\|_{TV} \leq 1.
\end{equation*}
Clearly, the solution of Ex. \ref{ex:basic}, with a total variation of $\frac{3}{2}$, violates this extra constraint, so the algorithm should converge to another solution. Again, from the first relaxation on, the cost of the associated truncated moment problem is $\frac{1}{8}$. It is also plain to see that $M_d \left( y^\mu, x  \right)$ is the truncated moment matrix of a Dirac located at $x=\frac{1}{4}$, hinting a trajectory $x(t)=\frac{1}{4}$. On the control side, starting from the second relaxation, it also becomes evident that $M_d \left( y^\nu, t  \right)$ is the truncated moment matrix of the signed measure $-\frac{3}{4} \, \delta_0 + \frac{1}{4} \, \delta_2$, revealing the times and amplitudes of impulses compatible with admissible controls. This leads to the trajectory of Fig. \ref{fig:totalVariation}, which therefore is an optimal solution of the problem.

\begin{figure}
  \centering
  \includegraphics[width=0.7\textwidth]{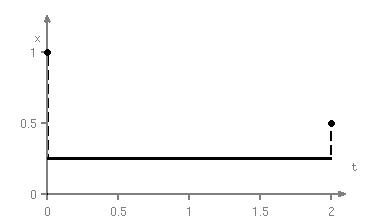}
  \caption{Trajectory for Ex. \ref{ex:totalVariation} }
  \label{fig:totalVariation}
\end{figure}

\end{example}

\begin{example}[Discrete control set with chattering]
\label{ex:DiscreteControl}
We take again Ex. \ref{ex:basic} with the additional constraint that the control $u(t)$ takes  its value in the set $ U=\left\lbrace \pm 1 \right\rbrace $, using the method explained in Section \ref{sec:discreteControl}. The solution to this problem is easy to infer: reach the turnpike $x(t)=0$ as quickly as possible by applying the negative control until $t=1$, then chatter with equal probability to remain on the turnpike until $t=\frac{3}{2}$, after which the positive control must be applied until $t=2$ (see Fig. \ref{fig:DiscreteControl}). This solution has an optimal cost of $\frac{3}{8} \approx 0.375$. Compare this value with those of Table \ref{tab:discreteControl}, which presents the evolution of the criterion with respect to the relaxation order of the truncated problem. After the fourth relaxation, the marginal w.r.t. $x$ of the control measure corresponding to the control $u(t)=+1$ closely approaches the positive measure $\frac{1}{2} \, I_{[1,\frac{3}{2}]}(dx) +  I_{[\frac{3}{2},2]}(dx)$ while the marginal w.r.t. $x$ of the control measure of $u(t)=-1$ 
converges to $ I_{[0,1]}(dx) + \frac{1}{2} \, I_{[1,\frac{3}{2}]}(dx) $, as expected.

\begin{figure}
  \centering
  \includegraphics[width=0.7\textwidth]{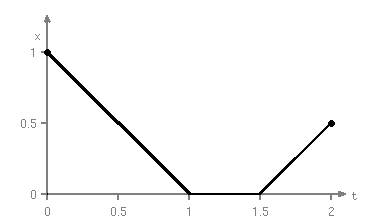}
  \caption{Trajectory for Ex. \ref{ex:DiscreteControl} }
  \label{fig:DiscreteControl}
\end{figure}

\begin{table}
  \renewcommand{\arraystretch}{1.3}
  \caption{Criterion as a function of LMI relaxation order for Ex. \ref{ex:DiscreteControl}}
  \label{tab:discreteControl}
  \centering
  \begin{tabular}{c||c|c|c|c}
$d$ & 1 & 2 & 3 & 4 \\     \hline
$V_M^d$ & 0.000 & 0.288  & 0.368 & 0.372 \\
  \end{tabular}
\end{table}

\end{example}

\begin{example}[Infeasible problem]
\label{ex:infeasible}
If the problem is infeasible, it \emph{may} be detected by the infeasibility of one of the LMI relaxations. Take Ex. \ref{ex:basic} with the additional total variation constraint $ \int_0^2 \! |u(t)| \, \dif t \leq \frac{1}{4}$ that puts the end point out of reach from the starting point. Indeed, at the first relaxation, the LMI problem is flagged as infeasible with
a Farkas dual vector, providing a certificate of infeasibility of the original problem.
\end{example}

\begin{example}[Unbounded problem]
\label{ex:unbounded}
If the problem is unbounded, it \emph{will} be detected at the first LMI relaxation.
Consider the problem of maximizing the total variation of a linear control problem:
\begin{equation}
  \sup_{u(t)} \int_0^1 \! |u(t)| \: \dif t
\end{equation}
such that
\begin{align*}
    & \dot{x}(t) = u(t) \\
    & x(0) = 0, \quad x(1) = 0 \\
    & x^2(t) \leq 1.
\end{align*}
As expected, the LMI problem from the first relaxation on is flagged as unbounded because its dual is flagged as infeasible.

\end{example}

\begin{example}[Bang-bang control of the Vanderpol equation]
\label{ex:VDP}

Consider the following time-optimal problem of the Vanderpol equation:

\begin{equation*}
  \inf_{u(t) \in U} T
\end{equation*}
such that
\begin{align*}
    & \dot{x}_1(t) = x_2(t) \\
    & \dot{x}_2(t) = -x_1(t)-(x_1^2(t)-1) \, x_2(t) + u(t) \\
    & x(0) = \begin{bmatrix}-0.4 & -0.6\end{bmatrix}', \quad x(T) = \begin{bmatrix}0.6 & 0.4 \end{bmatrix}'\\
    & U = \left\{ \pm 1 \right\}.
\end{align*}

In \cite{Simakov}, this problem is solved by applying a gradient-based optimization technique on a parametrization of admissible trajectories, with a minimum time of $2.14$. However, this method can only prove the local optimality of solutions. Applying our method, we obtain a cost of $2.15$ at the fifth relaxation, providing a (numerical) certificate of global optimality for that local solution.

\end{example}

\section{The fuel-optimal linear impulsive guidance rendezvous problem}

In this section, the proposed approach is applied to the far-range rendezvous in a linearised gravitational field. This problem is defined as a fixed-time minimum-fuel impulsive orbital transfer between two known circular orbits. Under Keplerian assumptions and for a circular rendezvous, the complete rendezvous problem may be decoupled between the
out-of-plane rendezvous problem for which an analytical solution may be found \cite{Carter91a} and the coplanar problem. Therefore, only coplanar circular rendezvous problems based on the Hill-Clohessy-Wiltshire equations and associated transition matrix \cite{Clohessy60} are considered for numerical illustration of
the proposed results. The general framework of the minimum-fuel fixed-time coplanar rendezvous problem in a linear setting is recalled in \cite{Carter91a} and \cite{Arzelier11} where an indirect method based on primer vector theory is proposed. Considering the necessity of easy-to-implement numerical solution for on-board guidance algorithms, direct methods based on linear programming (LP) problem may be used as in \cite{Mueller08}. For an \textit{a priori} fixed number of impulsive manoeuvres and using a classical transcription method \cite{Betts01} \cite{Louembet}, the genuine infinite-dimensional problem may be converted into a finite-dimensional approximation given by the following LP problem:
\begin{equation}
\begin{array}{lll}\label{probleme-RdV-N-fixe-anomalie}
V_{LP}= & \displaystyle\min_{u} ~\displaystyle\sum_{i=1}^{N}\|u_{\theta_i}\|_1\\
\mathrm{s.t.}    &   x(\theta_f) =  \Phi(\theta_f,\theta_1) x(\theta_1) + \displaystyle\sum_{i=1}^{N}\Phi(\theta_f,\theta_{i})Bu_{\theta_i}\\
&   x(\theta_1)=x_0,~x(\theta_f)=x_f
\end{array}
\end{equation}
where $\Phi$ is the Hill-Clohessy-Wiltshire transition matrix, $B=\left[\begin{array}{cc}0_{2\times 2} & 1_{2}\end{array}\right ]'$ and $u_{\theta_i}$ is the vector of velocity increments at $\theta_i$ in the local vertical local horizontal (LVLH) frame \cite{Arzelier11}. Time has been changed to the true anomaly $\theta$ for the independent variable as is usual in the literature \cite{Carter91a},
and it ranges in the interval $[\theta_1,\,\theta_f]$. Note that this formulation implies that only the impulsive solution of the general linear rendezvous problem may be obtained for a fixed number of velocity increments.

To be consistent with our previous notations we let $t=\theta$, $\theta_0=0$ and $\theta_f=T$.
Our impulsive optimal control problem (\ref{eq:OCP}) writes
\begin{equation*}
\begin{array}{lll}
V_M = & \displaystyle \inf_{w(t)} \int_0^T |dw_1|(t) + |dw_2|(t) \\
\mathrm{s.t.}    &  dx = \begin{bmatrix} 0 & 0 & 1 & 0 \\ 0 & 0 & 0 & 1 \\ 0 & 0 & 0 & 2 \\ 0 & 3 & -2 & 0\end{bmatrix} x(t) dt +
      \begin{bmatrix} 0 & 0 \\ 0 & 0 \\ 1 & 0 \\ 0 & 1 \end{bmatrix} dw(t) \\
    & x(0) = x_0, \quad x(T) = x_f
\end{array}
\end{equation*}
where state components model positions $(X,Z)=(x_1,x_2)$ in the orbital plane, and their respective velocities $(\dot{X},\dot{Z})=(x_3,x_4)$.
It clearly encompasses formulation (\ref{probleme-RdV-N-fixe-anomalie}) since it allows to consider continuous or impulsive thrusters as well. In both cases, the fuel consumption is measured by the one-norm of vector function $[\theta_1,\theta_f] \rightarrow \|dw\|_1$ \cite{Ross06} whereas the two-norm of this vector is used in general in the literature, see \cite{Carter91a}, \cite{Arzelier11} and references therein.

For the sake of comparison between these two approaches, two academic examples taken from \cite{Carter91a} are presented.

\begin{example}[In-plane rendezvous 1]
\label{ex:RDV_carter2}
Consider the first case presented in \cite{Carter91a}. It consists of a coplanar circle-to-circle rendezvous
with zero eccentricity. The
rendezvous manoeuvre must be completed in one orbital period with boundary conditions
$x_0=[ \begin{array}{cccc} 1 &  0& 0 & 0\end{array}]'$ and $x_f = [ \begin{array}{cccc} 0 & 0 & 0 & 0 \end{array}]'$.
This type of rendezvous is usually difficult to handle by numerical methods because of its singularity due to the high number of symmetries involved.

With a grid of $N=50$ points, the LP algorithm gives a two-impulse solution at the initial and final times of the rendezvous without interior impulse nor initial or final coasting period. The optimal impulses are both horizontal and opposite $u_0=-u_{2\pi}=\left[\begin{array}{cc}0.05305 & 0\end{array}\right ]'$. The fuel cost is given by $V_{LP}=0.1061$. 
The LMI method has no difficulty to recover the optimal solution given by the LP algorithm. A cost of $0.1061$ is obtained for each relaxation. It is then easy to extract from the matrices that the optimal solution for the first control consists of two symmetric impulses of magnitude $0.0531$ at the initial and final times, while the second control is identically $0$. The optimal trajectory in the orbital plane is depicted in Figure \ref{fig:RDV_Carter2_trajectoire} where + indicates the $50$ points of discretization.
\begin{figure}[H]
  \centering
  \includegraphics[width=0.8\textwidth]{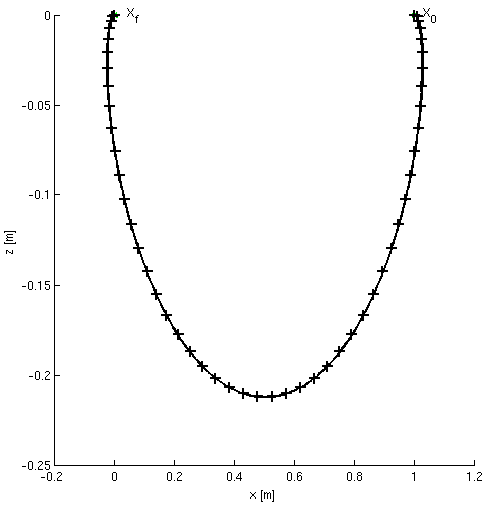}
  \caption{Trajectory in the orbital plane $(X,Z)$ in LVLH: Case 1 of \cite{Carter91a}}
  \label{fig:RDV_Carter2_trajectoire}
\end{figure}
Figure \ref{fig:RDV_Carter2_position} shows position, velocity and impulses history versus true anomaly.
\begin{figure}[H]
  \centering
  \includegraphics[width=0.7\textwidth]{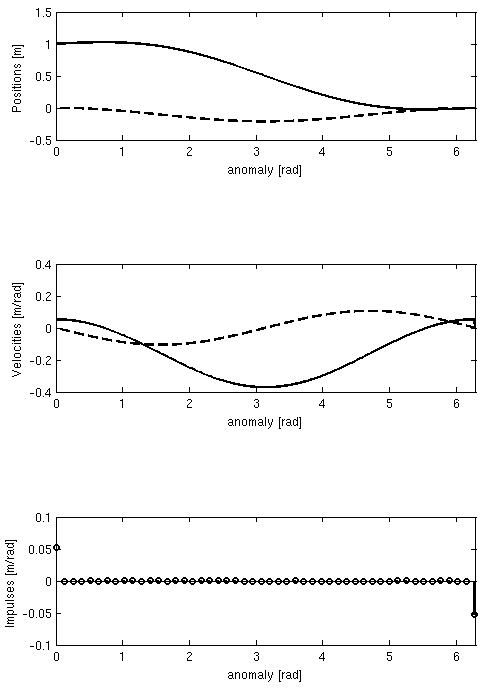}
  \caption{Positions ($X$ solid, $Z$ dashed), velocities ($\dot{X}$ solid, $\dot{Z}$ dashed) and impulses (on $X$ axis): Case 1 of \cite{Carter91a}}
  \label{fig:RDV_Carter2_position}
\end{figure}
\end{example}

\begin{example}[In-plane rendez-vous 2]
\label{ex:RDV_carter4}
As a second example, the third case of \cite{Carter91a} is revisited. The rendezvous is nearly identical to the previous one except for the final condition that imposes to reach the target with relative velocity of $0.427$ in the $Z$ direction, namely $x_0=[ \begin{array}{cccc} 1 &  0& 0 & 0
\end{array}]'$ and $x_f=[ \begin{array}{cccc} 0 & 0 & 0 & 0.427\end{array}]'$.

Again, a grid of $N=50$ points is used when running the LP algorithm. It converges to a four-impulse trajectory depicted in Fig. \ref{fig:RDV_Carter4_trajectoire}. The numerical results are summarized in Table \ref{tab:RDV_carter4}. 
\begin{figure}[H]
  \centering
  \includegraphics[width=0.8\textwidth]{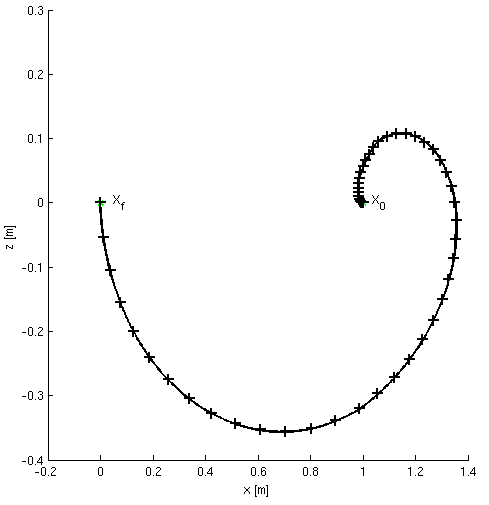}
  \caption{Trajectory in the orbital plane $(X,Z)$ in LVLH: Case 3 of \cite{Carter91a}}
  \label{fig:RDV_Carter4_trajectoire}
\end{figure}
\begin{table}[H]
  \renewcommand{\arraystretch}{1.4}
  \caption{Impulse times and amplitudes for Ex. \ref{ex:RDV_carter4}}
  \label{tab:RDV_carter4_sol}
  \centering
  \begin{tabular}{c|c|c||c|c|c}
    \hline
    \multicolumn{3}{c||}{LMI method} & \multicolumn{3}{c}{LP method} \\
    \hline
    \bfseries $\theta_i$ & \bfseries $(u_{\theta_i})_1$ & \bfseries $(u_{\theta_i})_2$ & \bfseries $\theta_i$ & \bfseries $(u_{\theta_i})_1$& \bfseries $(u_{\theta_i})_2$\\
    \hline\hline
    0 & -0.0386 & 0 & 0 & -0.0392& 0\\
    1.791 & +0.109 & 0 & 1.795 & +0.109 & 0\\
    4.495 & -0.109 & 0 & 4.488 & -0.109 & 0\\
    6.283 & +0.0389 & 0 & 6.283 & +0.0392 & 0\\
    \hline
  \end{tabular}
\end{table}
Using our algorithm, we reached the same criterion (within numerical tolerance) after the fourth relaxation (see Tab. \ref{tab:RDV_carter4}). As usual, the controls can be inferred from the moment matrix of the $\nu$ measures. Indeed, $\nu_1$ converges to the measure $\sum (u_{\theta_i})_1 \, \delta_{\theta_i}$ with impulse amplitudes $(u_{\theta_i})_1$ and anomaly $\theta_i$ taken from Table \ref{tab:RDV_carter4_sol}, while $\nu_2$ converges to an all zero measure. Not only does this result prove the global optimality of the conjectured solution within the class of all impulsive solutions no matter the number of impulses, but it also shows that it is optimal over all measure thrust solutions.
\begin{table}[H]
  \renewcommand{\arraystretch}{1.3}
  \caption{Criterion as a function of LMI relaxation order for Ex. \ref{ex:RDV_carter4}}
  \label{tab:RDV_carter4}
  \centering
  \begin{tabular}{c||c|c|c|c}
$d$ & 1 & 2 & 3 & 4 \\     \hline
$V_M^d$ & 0.0463  & 0.0680  & 0.2188 & 0.2972 \\
  \end{tabular}
\end{table}
Finally, position, velocity and impulses history are illustrated in Figure \ref{fig:RDV_Carter4_position}. Note the symmetry of the optimal four-impulse solution.
\begin{figure}[H]
  \centering
  \includegraphics[width=0.7\textwidth]{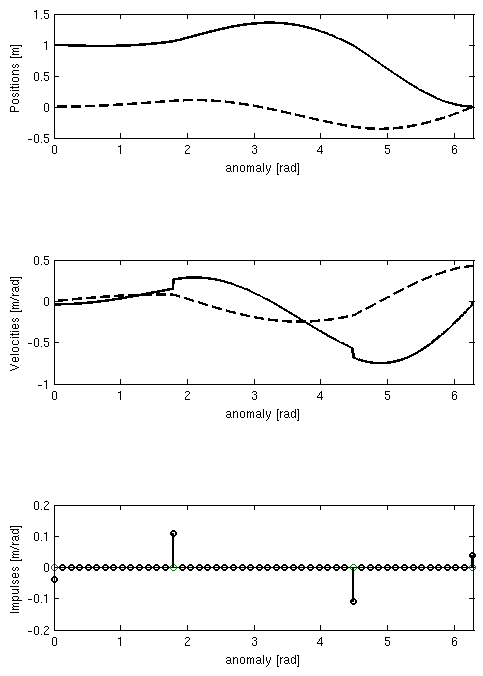}
  \caption{Positions ($X$ solid, $Z$ dashed), velocities ($\dot{X}$ solid, $\dot{Z}$ dashed) and impulses (on $X$ axis): Case 3 of \cite{Carter91a}}
  \label{fig:RDV_Carter4_position}
\end{figure}
\end{example}

\section{Conclusion}

The focus of this work is on actual computation of optimal impulsive controls
for systems described by ordinary differential equations with polynomial
dynamics and polynomial (semialgebraic) constraints on the state. State
trajectory and controls are measures which are linearly constrained,
resulting in an infinite-dimensional linear programming (LP) problem
consistent with the formalism of our GloptiPoly software \cite{GloptiPoly}.
This LP problem on measures can then be solved numerically via a
hierarchy of linear matrix inequality (LMI) relaxations, for which
off-the-shelf semi-definite programming (SDP) solvers can be used.
The optimal impulse sequence can then be retrieved by simple
linear algebra, and global optimality can be verified
by a posteriori simulation or comparison with suboptimal
control sequences computed by alternative techniques.

For space rendezvous, our technique can be readily adapted to cope with state (e.g. obstacle
avoidance) constraints, as soon as they are basic semialgebraic.
Other criteria than the total variation can also
be handled. Smoother solutions can be expected, maybe consisting
of a mix of absolutely continuous and singular controls, including
impulsive controls.

\end{document}